\documentclass[12pt]{article}
\usepackage{amsfonts}
\usepackage{latexsym, amssymb, amsmath, amscd, amsthm, amsfonts, mathrsfs, psfrag, ifpdf}
\usepackage{graphics}
\usepackage{graphicx}
\usepackage{color}

\newtheorem{theorem}{Theorem}[section]
\numberwithin{equation}{section}
\newtheorem{definition}[theorem]{Definition}
\newtheorem{lemma}[theorem]{Lemma}

\def\qed{\hfill \rule{4pt}{7pt}}

\def\pf{\noindent {\it Proof.} }

\begin{document}
\begin{center}
{\large
  {\bf Permutation Tableaux and

  the Dashed Permutation Pattern 32--1}}

\vskip 6mm {\small William Y.C. Chen$^1$  and Lewis H. Liu$^2$ \\[%
2mm] Center for Combinatorics, LPMC-TJKLC\\
Nankai University, Tianjin 300071,
P.R. China \\[3mm]
$^1$chen@nankai.edu.cn,  $^2$lewis@cfc.nankai.edu.cn \\[0pt%
] }
\end{center}


\begin{abstract}
We give a solution to a problem posed by Corteel and Nadeau concerning permutation tableaux of length $n$ and the number of occurrences of the dashed pattern 32--1 in permutations on $[n]$. We introduce the  inversion number
of a permutation tableau. For a permutation tableau $T$ and the permutation $\pi$ obtained from $T$ by the bijection of Corteel and Nadeau, we show that the inversion number of $T$ equals the number of occurrences of the dashed pattern 32--1 in the reverse complement of $\pi$. We also show that permutation tableaux without inversions coincide with L-Bell tableaux introduced by Corteel and Nadeau.
\end{abstract}




\section{Introduction}

Permutation tableaux were introduced by Steingr\'imsson and Williams \cite{Stei} in the study of totally positive Grassmannian cells \cite{lam, Postnikov, Williams}. They are closely related to the PASEP (partially asymmetric exclusion process) model in statistical physics \cite{Corteel1,new1, new2, ref4,vien2}. Permutation tableaux are also in one-to-one correspondence with alternative tableaux introduced by Viennot \cite{Viennot}.

A permutation tableau is defined by a Ferrers diagram possibly with empty rows such that the cells are filled with 0's and 1's subject to the following conditions:
\begin{enumerate}
\item Each column contains at least one 1.
\item There does not exist a $0$ with a $1$ above (in the same column)  and a $1$ to the left (in the same row).
\end{enumerate}

The length of a permutation tableau is defined as the number of rows plus the number of columns. A $0$ in a permutation tableau is said to be restricted if there is a $1$ above. Among the restricted 0's in a row, the rightmost one plays a special role, which is called a rightmost restricted 0. A row is said to be  unrestricted if it does not contain any restricted 0. A permutation tableau $T$ of length $n$ is labeled by the elements in $[n]=\{1, 2,\ldots,
n\}$ in increasing order from the top right corner to the bottom left corner. The set $[n]$ is referred to as  the label set of $T$. We use $(i,j)$ to denote the cell with row label $i$ and column label $j$.

For example, Figure 1.1 exhibits a permutation tableau of length 11 which contains an empty row. There are two rightmost restricted 0's at cells (5,9) and (8,10), and there are four unrestricted rows labeled by $1,2,7,$ and $11$.

\vspace*{4.5cm}
\begin{figure}[h]
\begin{center}
\hspace*{-3.5cm}
\begin{picture}(7,20)

\setlength{\unitlength}{0.7cm}

\put(0,7){\line(1,0){5}} \put(0,6){\line(1,0){5}} \put(0,5){\line(1,0){5}} \put(0,4){\line(1,0){3}}
\put(0,3){\line(1,0){2}} \put(0,2){\line(1,0){2}}

\put(0,7){\line(0,-1){6}} \put(1,7){\line(0,-1){5}} \put(2,7){\line(0,-1){5}} \put(3,7){\line(0,-1){3}}
\put(4,7){\line(0,-1){2}} \put(5,7){\line(0,-1){2}}

\put(0.5,6.5){\makebox(0,0){0}} \put(0.5,5.5){\makebox(0,0){0}}
\put(0.5,4.5){\makebox(0,0){0}} \put(0.5,3.5){\makebox(0,0){1}}
\put(0.5,2.5){\makebox(0,0){0}}

\put(1.5,6.5){\makebox(0,0){1}} \put(1.5,5.5){\makebox(0,0){1}}
\put(1.5,4.5){\makebox(0,0){0}} \put(1.5,3.5){\makebox(0,0){1}}
\put(1.5,2.5){\makebox(0,0){1}}

\put(2.5,6.5){\makebox(0,0){0}} \put(2.5,5.5){\makebox(0,0){0}}
\put(2.5,4.5){\makebox(0,0){1}}

\put(3.5,6.5){\makebox(0,0){0}} \put(3.5,5.5){\makebox(0,0){1}}

\put(4.5,6.5){\makebox(0,0){0}} \put(4.5,5.5){\makebox(0,0){1}}

\put(5.2,6.4){\makebox(0,0){\footnotesize 1}}
\put(5.2,5.4){\makebox(0,0){\footnotesize  2}}
\put(4.5,4.8){\makebox(0,0){\footnotesize  3}}
\put(3.5,4.8){\makebox(0,0){\footnotesize  4}}
\put(3.2,4.4){\makebox(0,0){\footnotesize  5}}
\put(2.5,3.8){\makebox(0,0){\footnotesize  6}}
\put(2.2,3.4){\makebox(0,0){\footnotesize  7}}
\put(2.2,2.4){\makebox(0,0){\footnotesize  8}}
\put(1.5,1.8){\makebox(0,0){\footnotesize  9}}
\put(0.5,1.8){\makebox(0,0){\footnotesize  10}}
\put(0.3,1.2){\makebox(0,0){\footnotesize  11}}
\put(2.5,0){\makebox(0,0){\mbox{ Figure 1.1: Permutation tableau }}}
\end{picture}
\end{center}
\end{figure}

It is known that the number of permutation tableaux of length $n$ is $n!$. There are several bijections between permutation tableaux and permutations, see Corteel and Nadeau \cite{Corteel}, Steingr\'imsson and Williams \cite{Stei}. The second bijection in \cite{Stei} connects the number of 0's in a permutation tableau to the total number of occurrences of the dashed patterns 31--2, 21--3 and 3--21. This bijection also yields a relationship between the number of 1's  and the number of occurrences of the dashed pattern 2--31 in a permutation, as well as a relationship between the number of cells in the Ferrers diagram and  the number of occurrences of dashed patterns 1--32 and 32--1 in a permutation. In answer to a question of Steingr\'imsson and Williams \cite{Stei}, Burstein \cite{Burstein} found a classification of zeros in permutation tableaux and its connection to the total number of occurrences of dashed patterns 31--2 and 21--3 , and the number of occurrences of the dashed pattern 3--21.

On the other hand, the second bijection of Corteel and Nadeau \cite{Corteel} implies that the number of non topmost 1's in a permutation tableau equals the number of occurrences of the dashed pattern 31--2 in the corresponding permutation. They raised the problem of finding a statistic on permutation tableaux that has the same
distribution as the number of occurrences of the dashed pattern 32-1 in permutations.

Let us recall the definition of  dashed permutation patterns introduced by Babson and Steingr\'imsson \cite{Babson}. A dashed pattern is a permutation of $[k]$ for $k\leq n$ containing dashes indicating that the entries in a permutation  of $[n]$ need not occur consecutively. In this notation, a permutation pattern $\sigma = \sigma_1 \sigma_2 \ldots \sigma_k$ in the usual sense may be rewritten as $\sigma = \sigma_1$--$\sigma_2$--$\ldots$--$\sigma_k$. For example, we say that a permutation $\pi$ on $[n]$ avoids a dashed pattern 32--1 if there are no subscripts $i<k$ such that $\pi_{i-1}>\pi_{i}>\pi_k$. Claesson and Mansour \cite{Claesson} found explicit formulas for the number of permutations containing exactly $i$ occurrences of a dashed pattern $\sigma$ of length 3 for $i=1,2,3$.

The main idea of this paper is to introduce the inversion number of a permutation tableau (see Definition \ref{defi}). We show that the inversion number of a permutation tableau of length $n$ has the same distribution as the number of occurrences of dashed pattern 32-1 in a permutation of $[n]$. To be more specific, for a permutation tableau $T$ and the permutation $\pi$ obtained from $T$ by the first bijection of Corteel and Nadeau, we prove that the inversion number of $T$ equals the number of occurrences of the dashed pattern 32--1 in the reverse complement  of $\pi$. This gives a solution to the problem proposed Corteel and
Nadeau \cite{Corteel}.

This paper is ended with a connection between permutation tableaux without inversions  and L-Bell tableaux introduced by Corteel and Nadeau in \cite{Corteel}. We show that a permutation tableau has no  inversions  if and only if it is an L-Bell tableau.


\section{The inversion number of a permutation tableau }

In this section, we define the inversion number of a permutation
tableau (see Definition \ref{defi}). We show that the inversion
number of permutation tableau $T$ equals the number of occurrences of
the dashed pattern 32--1 in the reverse complement
of the permutation $\pi$ corresponding to $T$ under the bijection $\xi$ of Corteel and Nadeau.

Let ${\pi}=\pi_1\pi_2\ldots \pi_n$. Denote by $f_{\sigma}(\pi)$ the number of occurrences of dashed pattern $\sigma$ in $\pi$. The reverse complement of $\pi$ is denoted by
\[ \bar{\pi}= (n+1-\pi_n, \ldots n+1-\pi_{2}, n+1-\pi_1),\]
 here a permutation is also written in the
form of a vector.

The  main result of this paper is  stated as follows. Throughout this paper $\xi$ denotes the
bijection of Corteel and Nadeau. For completeness, a brief description of $\xi$ will be given later.

\begin{theorem}\label{main}
Let $T$ be a permutation tableau. Let $\mbox{inv}(T)$  be the number of inversions of $T$. Then we have
\begin{equation} \label{if}
\mbox{inv}(T)=f_{32\text{--}1}(\bar{\pi}).
\end{equation}
\end{theorem}

Since an occurrence  of the dashed pattern 32--1 in
$\bar{\pi}$ corresponds to an occurrence  of the dashed pattern 3-21 in ${\pi}$,   relation (\ref{if}) can be restated as
\begin{equation}\label{if-2}
\mbox{inv}(T)=f_{3\text{--}21}( {\pi}).
\end{equation}

Corteel and Nadeau \cite{Corteel} have shown that a permutation tableau is uniquely determined by the
topmost 1's and the rightmost restricted 0's. For a permutation tableau $T$, the alternative representation of $T$ is defined by Corteel and Kim \cite{Corteel3} as the diagram
obtained from $T$ by replacing the topmost 1's with $\uparrow$'s, replacing the rightmost restricted 0's
with $\leftarrow$'s and leaving the remaining cells blank. In this paper, we use black dots and white dots to represent the topmost 1's and rightmost restricted 0's in an alternative representation, as illustrated in Figure 2.1.

Note that for each dot in an alternative representation of a permutation tableau, there is an unique path, called an alternating path, to a black dot on an  unrestricted row. More precisely, an alternating paths is defined as follows. For a white dot we can find a black dot at the top as the next dot. For a black dot which is not on an unrestricted row, there is a white dot to the left as the next dot. Figure 2.1 shows two alternating paths.

For an alternative representation of a permutation tableau, we may use the row and column labels to represent an alternating path. It is evident that a black dot is
determined by a column label and a white dot is determined by a row label. Hence
an alternating path can be represented by an alternating sequence of row and column labels ending with a column label of a black dot on an unrestricted row.

For example, for the black dot in cell $(5,6)$, the alternating path can be expressed as  $(6, 5, 12).$ For the white dot in cell $(7,10)$, the corresponding alternating path takes the form $(7,10,4,11).$

\vspace*{4.5cm}
\begin{figure}[h]
\begin{center}
\hspace*{-10.6cm}
\begin{picture}(7,20) \setlength{\unitlength}{0.7cm}

\put(0,7){\line(1,0){6}} \put(0,6){\line(1,0){6}} \put(0,5){\line(1,0){5}} \put(0,4){\line(1,0){5}}
\put(0,3){\line(1,0){5}} \put(0,2){\line(1,0){4}}
\put(0,1){\line(1,0){3}}

\put(0,7){\line(0,-1){6}} \put(1,7){\line(0,-1){6}} \put(2,7){\line(0,-1){6}} \put(3,7){\line(0,-1){6}}
\put(4,7){\line(0,-1){5}} \put(5,7){\line(0,-1){4}}
\put(6,7){\line(0,-1){1}}

\put(0.5,6.5){\makebox(0,0){1}}
\multiput(0.5,5.5)(0,-1){5}{\makebox(0,0){0}}

\put(1.5,6.5){\makebox(0,0){1}} \put(1.5,5.5){\makebox(0,0){1}}
\put(1.5,4.5){\makebox(0,0){0}} \put(1.5,3.5){\makebox(0,0){1}}
\multiput(1.5,2.5)(0,-1){2}{\makebox(0,0){0}}

\multiput(2.5,6.5)(0,-1){2}{\makebox(0,0){0}}
\multiput(2.5,4.5)(0,-1){2}{\makebox(0,0){1}}
\multiput(2.5,2.5)(0,-1){2}{\makebox(0,0){0}}

\multiput(3.5,6.5)(0,-1){2}{\makebox(0,0){0}}
\multiput(3.5,4.5)(0,-1){3}{\makebox(0,0){1}}

\multiput(4.5,6.5)(0,-1){3}{\makebox(0,0){0}}
\put(4.5,3.5){\makebox(0,0){1}} \put(5.5,6.5){\makebox(0,0){1}}

\put(6.2,6.4){\makebox(0,0){\footnotesize 1}}
\put(5.5,5.75){\makebox(0,0){\footnotesize  2}}
\put(5.2,5.4){\makebox(0,0){\footnotesize  3}}
\put(5.2,4.4){\makebox(0,0){\footnotesize  4}}
\put(5.2,3.4){\makebox(0,0){\footnotesize  5}}
\put(4.5,2.75){\makebox(0,0){\footnotesize  6}}
\put(4.2,2.4){\makebox(0,0){\footnotesize  7}}
\put(3.5,1.75){\makebox(0,0){\footnotesize  8}}
\put(3.2,1.4){\makebox(0,0){\footnotesize  9}}
\put(2.5,0.75){\makebox(0,0){\footnotesize  10}}
\put(1.5,0.75){\makebox(0,0){\footnotesize 11}}
\put(0.5,0.75){\makebox(0,0){\footnotesize  12}}\hspace*{7cm}
\put(0,7){\line(1,0){6}} \put(0,6){\line(1,0){6}} \put(0,5){\line(1,0){5}} \put(0,4){\line(1,0){5}}
\put(0,3){\line(1,0){5}} \put(0,2){\line(1,0){4}}
\put(0,1){\line(1,0){3}}

\put(0,7){\line(0,-1){6}} \put(1,7){\line(0,-1){6}} \put(2,7){\line(0,-1){6}} \put(3,7){\line(0,-1){6}}
\put(4,7){\line(0,-1){5}} \put(5,7){\line(0,-1){4}}
\put(6,7){\line(0,-1){1}}

\put(0.5,6.5){\makebox(0,0){$\bullet$}}
\multiput(0.5,5.5)(0,-2){2}{\makebox(0,0){$\circ$}}

\put(1.5,6.5){\makebox(0,0){$\bullet$}}
\put(1.5,4.5){\makebox(0,0){$\circ$}}

\multiput(2.5,4.5)(0,-1){1}{\makebox(0,0){$\bullet$}}
\multiput(2.5,2.5)(0,-1){2}{\makebox(0,0){$\circ$}}

\multiput(3.5,4.5)(0,-1){1}{\makebox(0,0){$\bullet$}}

\put(4.5,3.5){\makebox(0,0){$\bullet$}}

\put(5.5,6.5){\makebox(0,0){$\bullet$}}

 \put(.5,3.5){\line(1,0){4}}
\put(1.5,4.5){\line(1,0){1}}
\put(0.5,6.5){\line(0,-1){3}} \put(1.5,6.5){\line(0,-1){2}}
\put(2.5,4.5){\line(0,-1){2}}

\put(.5,6.2){\makebox(0,0){$\uparrow$}}
\put(0.78,3.465){\makebox(0,0){$\leftarrow$}}
\put(2.5,4.2){\makebox(0,0){$\uparrow$}}
\put(1.78,4.465){\makebox(0,0){$\leftarrow$}}
\put(1.5,6.2){\makebox(0,0){$\uparrow$}}

\put(6.2,6.4){\makebox(0,0){\footnotesize 1}}
\put(5.5,5.75){\makebox(0,0){\footnotesize  2}}
\put(5.2,5.4){\makebox(0,0){\footnotesize  3}}
\put(5.2,4.4){\makebox(0,0){\footnotesize  4}}
\put(5.2,3.4){\makebox(0,0){\footnotesize  5}}
\put(4.5,2.75){\makebox(0,0){\footnotesize  6}}
\put(4.2,2.4){\makebox(0,0){\footnotesize  7}}
\put(3.5,1.75){\makebox(0,0){\footnotesize  8}}
\put(3.2,1.4){\makebox(0,0){\footnotesize  9}}
\put(2.5,0.75){\makebox(0,0){\footnotesize  10}}
\put(1.5,0.75){\makebox(0,0){\footnotesize 11}}
\put(0.5,0.75){\makebox(0,0){\footnotesize  12}}

\put(-2.8,-.5){\makebox(0,0){\mbox{ Figure 2.1: Alternative representation of a
permutation tableau}}}
\end{picture}
\end{center}
\end{figure}

For two alternating paths $P$ and $Q$ of $T$, we say that $P$ contains $Q$ if $Q$ is a segment of $P$. Assume that $P$ and $Q$ does not contain each other. We proceed to define the order relation on $P$ and $Q$. Clearly, if $P$ and $Q$ intersect
at some point, then they will share the same
ending segment after this point. If this is the case, we will remove the common dots of $P$ and $Q$, and then consider the resulting
alternating paths $P'$ and $Q'$. We say that   $P>Q$ ($P<Q$) if one of the following four conditions holds:
\begin{enumerate}
\item  The paths $P$ and $Q$ do not intersect, and the  ending  dot of $P$ is below (above) the ending dot of  $Q$.

\item  The paths $P$ and $Q$ do not intersect, and the ending  dot of $P$ is to the right (to the left) of the ending dot of $Q$.

\item The paths $P$ and $Q$ intersect, and the ending dot of $P'$ is below (above) the ending dot of $Q'$.

\item The paths $P$ and $Q$ intersect, and the ending point of $P'$ is to the right (to the left) of the ending dot of $Q'$.
\end{enumerate}

For any two alternating paths $P$ and $Q$ that do  not contain each other, it can be seen that either $P>Q$ or $P<Q$ holds.
Using this order, we can define the inversion number of an alternative representation $T$ of a
permutation tableau, which is also considered as the inversion number of the original permutation tableau.

\begin{definition}\label{defi}
Suppose that $j$ is a column label of $T$ and $P_j$ is the alternating path  starting with the black dot
with column label $j$. Let $k$ be a label of $T$ with $j<k$ such that $k$ is not a label of a dot on $P_j$. In other words, $P_k$ is not contained
in $P_j$, where $P_k$ denotes the alternating path with the dot with label $k$. We say that the pair of labels $(j,k)$ is an inversion of $T$ if $P_j>P_k$. The total number of inversions of $T$ is denoted by
$\rm{inv}(T)$.
\end{definition}

Analogous to the inversion code of a permutation, for a column label $j$ we define $w_j(T)$ as the number
of inversions of $T$ that are of the form $(j,k)$. Hence
$$\mbox{inv}(T)=\sum_{j\in C(T)}w_j(T),$$
where $C(T)$ is the set of column labels of $T$.

For example, Figure 2.2 gives two permutation tableaux in the form of their alternative representations. For the  alternative representation $T$ on the left, we have $C(T)=\{2,3\}$. Since $P_2>P_3$, we see that  $w_2(T)=1$, $w_3(T)=0$, and   $\mbox{inv}(T)=1.$ For the alternative representation $T'$ on the right, we have $C(T')=\{3,5\}$. Since  $P_3>P_4$ and $P_3>P_5$, we find $w_3(T')=2$, $w_5(T')=0$, and $\mbox{inv}(T')=2.$

\vspace*{3.5cm}
\begin{figure}[h]
\vspace*{-1.5cm} \hspace*{3.5cm}
\begin{picture}(4,5)
\setlength{\unitlength}{0.7cm}
\multiput(0,2)(0,-1){2}{\line(1,0){2}}
\multiput(0,2)(1,0){3}{\line(0,-1){1}}

\put(.5,1.5){\makebox(0,0){$\bullet$}}
\put(1.5,1.5){\makebox(0,0){$\bullet$}}

\put(2.3,1.5){\makebox(0,0){$1$}} \put(1.5,.7){\makebox(0,0){$2$}}
\put(0.5,.7){\makebox(0,0){$3$}} \hspace*{6cm}
\multiput(0,2)(0,-1){3}{\line(1,0){2}}
\multiput(0,-1)(0,-1){1}{\line(1,0){1}}
\multiput(0,2)(1,0){2}{\line(0,-1){3}}
\multiput(2,2)(1,0){1}{\line(0,-1){2}}

\put(.5,1.5){\makebox(0,0){$\bullet$}}
\put(.5,-.5){\makebox(0,0){$\circ$}}
\put(1.5,.5){\makebox(0,0){$\bullet$}}
\put(.5,1.2){\makebox(0,0){$\uparrow$}}

\put(.5,1.5){\line(0,-1){2}}

\put(2.3,1.5){\makebox(0,0){$1$}} \put(2.3,.5){\makebox(0,0){$2$}}
\put(1.7,-.35){\makebox(0,0){$3$}}
\put(1.3,-.68){\makebox(0,0){$4$}} \put(.6,-1.2){\makebox(0,0){$5$}}

\put(-3.5,-2){\makebox(0,0){\mbox{ Figure 3.2: Two examples}}}
\end{picture}

\end{figure}
\vspace*{1.5cm}

We now proceed to present a proof of Theorem 2.1. For completeness, we give a brief description of the bijection $\xi$ of Corteel and Nadeau from permutation tableaux to permutations.

Assume that $T$ is the alternative representation of a permutation tableau. Let $\xi(T)=\pi=\pi_1\pi_2\ldots\pi_n$. The bijection is a recursive procedure to construct $\pi$.
Starting with the sequence of the labels of free rows  in
increasing order. Then  successively insert the column labels of $T$. Let $j$ be    the maximum column label to be inserted.  If the cell
$(i,j)$ is filled with a black dot, then insert $j$ immediately to the
left of $i$. If column $j$ contains white dots in rows
$i_1,i_2,\ldots,i_k$, then insert $i_1,i_2,\ldots,i_k$ in increasing
order to the left of $j$. Repeating this process, we obtain a  permutation $\pi$.

For the sake of presentation, we give
 three lemmas that will be needed in the proof
 of theorem \ref{main}.
The first lemma was established by Corteel and Nadeau \cite{Corteel}.

\begin{lemma}\label{lem}
Let $\pi=\xi(T)$.  Then $\pi_{i}>\pi_{i+1}$  if
and only $\pi_i$ is a column label of $T$.
\end{lemma}

The next lemma states that the labels representing an alternating path of $T$ form a subsequence of  $\xi(T)$.

\begin{lemma}\label{lemm}
Let $P=p_1p_2\cdots p_r$ be an alternating path of $T$ starting with  a dot labeled by $p_1$ and ending with a black dot labeled by $p_r$. Then
$p_1p_2\cdots p_r$ is a subsequence of  $\xi(T)$.
\end{lemma}

\pf Assume that the alternating path $P$ ends with a black dot at cell $(i,p_r)$, where $i$ is an unrestricted row label.  Since the ending dot
 represents a topmost 1, by the construction of $\xi$, we see that  $p_r$ is inserted to the left of $i$. Note that the cell $(p_{r-1},p_r)$ is filled with a white dot representing a rightmost restricted 0, so $p_{r-2}$ is inserted to the left of $p_{r-1}$. Since the path $P$ is alternating with respect to black and white dots, we deduce that $p_{r-3},\cdots,p_2,p_1$ are inserted one by one such that $p_{i}$ is inserted to the left of $p_{i+1}$ for $i=1,2,\ldots,r-1$. It follows that $p_1p_2\cdots p_r$ is a subsequence of the permutation $\xi(T)$. This completes the proof. \qed

 Given two labels $i$ and $j$ of $T$, the following lemma shows that the relative order of $i$ and $j$ in $\xi(T)$ can be determined by
 the order of the alternating paths starting with the two dots labeled by $i$ and $j$.

\begin{lemma}\label{prop}
Let $P_{i}$ and $P_j$ be two alternating paths of $T$ starting with two dots labeled by $i$ and $j$ respectively. Then $i$ is to the left of $j$ in $\xi(T)$ if and only if $P_i<P_j$ or $P_j$ is contained in $P_i$.
\end{lemma}

\pf If $P_j$ is contained in $P_i$, by Lemma \ref{lemm} we see that $i$ is to the left of $j$. Otherwise, let $P_{i}=i_1i_2\ldots i_s$ and $P_{j}=j_1j_2\ldots j_t$, where  $i=i_1$ and $j=j_1$ respectively. Assume that $P_{i}$ and $P_{j}$ do not intersect. We wish to prove that $j_1$ is to the right of $i_s$ in $\xi(T)$, since  $i_1$ is to the left of $i_s$ by Lemma \ref{lemm}.

 If  $P_{i}$ and $P_{j}$ intersect at some dot, let $P_{i}'$ and $P_{j}'$ be the alternating paths by removing the common dots of $P_{i}$ and $P_{j}$. Suppose that the ending dot
  of $P_{i}'$ is labeled by $i_{s-m}$. Then the ending dot of $P_{j}'$ is labeled by $j_{t-m}$. In this case, we aim to show that $j_1$ is to the right of $i_{s-m}$ in $\xi(T)$, since $i_1$ is to the left of $i_{s-m}$ by Lemma \ref{lemm}.

 Consider the following four cases corresponding to the order on  $P_{i}$ and $P_{j}$. In the first two cases, $P_{i}$ and $P_{j}$ do not intersect, that is, $i_m\neq j_r$ for any $m=1,2,\ldots,s$ and $r=1,2,\ldots,t$. Suppose that the ending black dots of $P_i$ and $P_j$ are in cells $(r_{P_i},i_s)$ and $(r_{P_j},j_t)$ respectively, where $r_{P_i}$ and $r_{P_j}$ are the labels of unrestricted rows.

\noindent
Case 1. $r_{P_i}<r_{P_j}$. First of all, $r_{P_i}$ is to the left of $r_{P_j}$ in $\xi(T)$. Since both cells $(r_{P_i},i_s)$ and $(r_{P_j},j_t)$ are filled
with black dots, the element $i_s$ is inserted immediately to the left of $r_{P_i}$, while $j_t$ is inserted immediately to the left of $r_{P_j}$. This implies that  $j_t$ is to the right of $r_{P_i}$ in $\xi(T)$. Hence $j_t$ is to the right of $i_s$.

Since the path $P_j$ alternates with black  and white dots, the cell $(j_t,j_{t-1})$ is filled with a white dot. Thus $j_{t-1}$ is inserted to the left
of $j_t$ but to the right of $r_P$, that is, $j_{t-1}$ is to the right of $i_s$. Repeating the above procedure, we reach the conclusion that  the label $j_r$ is to the right of $i_s$ for $r=t,t-1, \ldots,1$. In particular,  $j_1$ is to the right of $i_s$, so that $j_1$ is to the right of $i_1$.

\noindent
Case 2.  $r_{P_i}=r_{P_j}$ and $j_t<i_s$. In the implementation of the  algorithm $\xi$, $i_s$ is inserted immediately to the left of $r_{P_i}$ and then $j_t$ is inserted immediately to the left of $r_{P_j}=r_{P_i}$. Hence $j_t$ is to the right of $i_s$. Inspecting the relative positions  of $i_s$ and
$j_r$ for $r<t$ like we have done in Case 1, we infer that $j_r$
is to the right of $i_s$ for $r=t,t-1,\ldots,1$.
So we arrive at the assertion that $j_1$ is to the right of $i_1$.

\noindent
Case 3. The ending point of $P_{j}'$ is below the ending point of $P_{i}'$, that is, $j_{t-m}>i_{s-m}$. In this case, both cells $(i_{s-m},i_{s-m+1})$
and $(j_{t-m},j_{t-m+1})$ are filled with white dots. To construct $\pi$ from $T$ according to $\xi$, both elements  $i_{s-m}$ and $j_{t-m}$ are inserted  to the left of the element $i_{s-m+1}=j_{t-m+1}$ in increasing order. Hence $i_{t-m}$ is to the left of $j_{t-m}$. Considering the relative positions of $i_{s-m}$ and $j_r$
for $r<t-m$ as in Case 1, we deduce that $j_r$
is to the right of $i_{s-m}$ for $1\leq r<t-m$. Therefore $j_1$ is to the right of $i_{s-m}$, and hence to the right of $i_1$.

\noindent
Case 4. The ending point of $P_{j}'$ is to the right of the ending point of $P_{i}'$, that is, $j_{t-m}<i_{s-m}$. Observe  that the element
$i_{s-m}$ is inserted immediately to the left of $i_{s-m+1}$. Then the element $j_{t-m}$ is inserted
immediately to the left of $j_{t-m+1}=i_{s-m+1}$. It follows that  $j_{t-m}$ is to the right of $i_{s-m}$. Using the same argument as in Case 1 for the elements  $i_{s-m}$ and $j_r$ with $r<t-m$, we conclude that $j_r$
is to the right of $i_{s-m}$ for $1\leq r<t-m$. Consequently, $j_1$ is to the right of $i_{s-m}$, and hence to the right of $i_1$.

 In summary, we see that if $P_{i}<P_{j}$, then $i$ is to the left of $j$ in $\xi(T)$. It remains to show that if $i$ is to the left of $j$ in $\xi(T)$, then we have $P_i<P_j$ or $P_j$ is contained in $P_i$. The proof is essentially the reverse procedure of the above argument, and is omitted.
\qed

\noindent
{\textit{Proof of Theorem \ref{main}.}}
Let $\xi(T)=\pi=\pi_1\pi_2\ldots\pi_n$. Combining Lemma \ref{lem} and Lemma \ref{prop}, we find that the subsequence
$\pi_{i}\pi_{j}\pi_{j+1}$ of $\pi$ is an occurrence  of the dashed
pattern 3--21  if and only if $(\pi_j,\pi_i)$ is an inversion of $T$. It follows that  $$\mbox{inv}(T)=f_{3\text{--}21}( {\pi}),$$
which is equivalent to the statement of the theorem.  \qed

  Figure 3.4 gives a permutation
tableau $T$ of length 10 and its alternative representation. For this example, we see that  $C(T)=\{5,6,7,9\}$,
 $w_5(T)=4$, $w_6(T)=3$, $w_8(T)=1$
 and $w_9(T)=0$. Hence we have $\mbox{inv}(T)=8$.
On the other hand,
$$\pi=\xi(T)=(9,2,7,8,1,6,5,3,4).$$
 This gives $\bar{\pi}=(6,5,7,4,9,2,3,8,1).$ It can be checked that the number of occurrences of the dashed  pattern
32--1 in $\bar{\pi}$ equals 8.

\vspace*{4.5cm}
\begin{figure}[h]
\begin{center}
\hspace*{-8.3cm}
\begin{picture}(7,20)

\setlength{\unitlength}{0.7cm}

\multiput(0,6)(0,-1){5}{\line(1,0){4}}
\put(0,1){\line(1,0){2}}

\multiput(0,6)(1,0){3}{\line(0,-1){5}} \multiput(3,6)(1,0){2}{\line(0,-1){4}}

\put(0.5,5.5){\makebox(0,0){1}} \multiput(0.5,4.5)(0,-3){2}{\makebox(0,0){0}} \multiput(0.5,3.5)(0,-1){2}{\makebox(0,0){1}}

\put(1.5,5.5){\makebox(0,0){1}} \put(1.5,4.5){\makebox(0,0){0}} \put(1.5,3.5){\makebox(0,0){1}} \put(1.5,2.5){\makebox(0,0){1}} \put(1.5,1.5){\makebox(0,0){0}}

\multiput(2.5,5.5)(0,-1){2}{\makebox(0,0){0}} \multiput(2.5,3.5)(0,-1){2}{\makebox(0,0){1}}

\multiput(3.5,5.5)(0,-1){2}{\makebox(0,0){0}} \multiput(3.5,3.5)(0,-1){2}{\makebox(0,0){1}}

\put(4.2,5.5){\makebox(0,0){\footnotesize 1}}    \put(4.2,4.5){\makebox(0,0){\footnotesize  2}}
\put(4.2,3.5){\makebox(0,0){\footnotesize  3}}   \put(4.2,2.5){\makebox(0,0){\footnotesize  4}}
\put(3.5,1.75){\makebox(0,0){\footnotesize  5}}   \put(2.5,1.75){\makebox(0,0){\footnotesize  6}}
\put(2.2,1.5){\makebox(0,0){\footnotesize  7}}   \put(3.5,1.75){\makebox(0,0){\footnotesize  8}}
\put(1.5,0.75){\makebox(0,0){\footnotesize 8}}  \put(0.5,0.75){\makebox(0,0){\footnotesize  9}}\hspace*{6cm}

\multiput(0,6)(0,-1){5}{\line(1,0){4}}
\put(0,1){\line(1,0){2}}

\multiput(0,6)(1,0){3}{\line(0,-1){5}} \multiput(3,6)(1,0){2}{\line(0,-1){4}}

\put(0.5,5.5){\makebox(0,0){$\bullet$}}

\put(1.5,5.5){\makebox(0,0){$\bullet$}} \put(1.5,4.5){\makebox(0,0){$\circ$}} \put(1.5,1.5){\makebox(0,0){$\circ$}}

\multiput(2.5,3.5)(0,-1){1}{\makebox(0,0){$\bullet$}}

\multiput(3.5,3.5)(0,-1){1}{\makebox(0,0){$\bullet$}}

\put(1.5,5.5){\line(0,-1){4}}
\put(1.5,5.2){\makebox(0,0){$\uparrow$}}
\put(4.2,5.5){\makebox(0,0){\footnotesize 1}}    \put(4.2,4.5){\makebox(0,0){\footnotesize  2}}
\put(4.2,3.5){\makebox(0,0){\footnotesize  3}}   \put(4.2,2.5){\makebox(0,0){\footnotesize  4}}
\put(3.5,1.75){\makebox(0,0){\footnotesize  5}}   \put(2.5,1.75){\makebox(0,0){\footnotesize  6}}
\put(2.2,1.5){\makebox(0,0){\footnotesize  7}}   \put(1.5,0.75){\makebox(0,0){\footnotesize 8}}
\put(0.5,0.75){\makebox(0,0){\footnotesize  9}}
\put(-2.8,-.5){\makebox(0,0){\mbox{ Figure 3.2: A permutation tableau and its alternative representation}}}
\end{picture}
\end{center}
\end{figure}

\section{Connection with L-Bell tableaux}

In this section, we show that a
permutation tableau has no inversions if and only if it is an L-Bell tableau as  introduced by Corteel and Nadeau \cite{Corteel}.
Recall that an L-Bell tableau is a permutation tableau such that any topmost 1 is also a left-most 1.

It has been shown by Claesson \cite{Claesson0} that the number of permutations of $[n]$ avoiding the dashed pattern 32--1 is given by the $n$-th Bell number $B_n$.  Together with Theorem \ref{main} we are led to the following correspondence.

\begin{theorem}\label{lem2}
The number of permutation tableaux $T$ of length $n$ such that
$\mbox{inv}(T)=0$ equals  $B_n$.
\end{theorem}

On the other hand, the following relation was
 proved by Corteel and Nadeau \cite{Corteel}.

\begin{theorem}\label{lem3}
The number of L-Bell permutation tableaux  of length $n$ equals $B_n$
\end{theorem}

By the definition of an inversion of a permutation tableau, it is straightforward to check that
an L-Bell tableau has no inversion. Combining
Theorem 3.2 and Theorem 3.3, we obtain the
following connection.

\begin{theorem}
Let $T$ be a permutation tableau. Then  $\mbox{inv}(T)=0$ if and only if $T$ is an L-Bell tableau.
\end{theorem}

Here we give a direct reasoning of the above theorem. Let $T$ be an alternative representation of a permutation tableau without inversions. It can be seen that the  permutation tableau corresponding to $T$ is an L-Bell tableau if and only  $T$ satisfies the following conditions:
\begin{enumerate}
\item Each row contains at most one black dot.

\item There is no empty cell that has a black dot above and a black dot immediately to the right.
\end{enumerate}
We wish to prove that if $\mbox{inv}(T)=0$, then $T$ satisfies the above conditions.

  Assume that there is a row containing two black dots, say, at  cells $(i,j)$ and $(i,k)$ with $j<k$. Then $(j,k)$ is an inversion of $T$, a contradiction. Thus condition (1) holds.

 For condition (2), assume to the contrary that there exists an empty cell $(i,j)$ such that there is a black dot above and a black dot immediately to the right. Without loss of generality, we may assume that $i$ is the minimal row label of an empty cell subject to the above assumption. Now we choose $j$ to be the maximal column label. Assume that the black dot in row $i$ and the black dot in column $j$ appear  at cells $(i,k)$ and $(t,j)$ respectively. Note that we have $t<i$. Since $j<k$ and $\mbox{inv}(T)=0$, we must have $P_{j}<P_{k}$. This implies that the unrestricted row containing the ending black dot of $P_{j}$ must be above row $t$. Assume that this ending black dot is in column $m$. Then there exists a white dot in row $t$. Suppose that it occurs at   cell $(t,s)$. Clearly, we have $s\leq m$. Moreover, we see that no cell $(x,y)$ for $x>t, y>s$ can  be filled with a white dot. On the other hand, no cell $(x,y)$ for $x<t,j<y<s$ can  be filled with a black dot.

By walking backwards from the ending black dot along the path $P_j$, we can find a black dot in column $s$  on the alternating path $P_j$ which is also on the alternating path $P_k$. So we deduce that the next white dot on $P_j$  is below the next white dot on $P_k$, that is, $P_j>P_k$. But this implies that  $(j,k)$ is an inversion of $T$, again a contradiction. Hence the proof is complete.\qed

\vspace{0.6cm} \noindent{\bf Acknowledgments.}  This work was
supported by the 973 Project, the PCSIRT Project of the Ministry of
Education, and the National Science Foundation of China.

\end{document}